\documentclass[letterpaper,10pt,conference]{ieeeconf}
\IEEEoverridecommandlockouts
\newtheorem{theorem}{Theorem}[section]
\newtheorem{lemma}[theorem]{Lemma}
\newtheorem{algorithm}[theorem]{Algorithm}

\newenvironment{remark}{
       \refstepcounter{theorem}\begin{trivlist}\item[]{\bfseries
       Remark \thetheorem\,}}
       {\end{trivlist}}

\usepackage[english]{babel}
\usepackage{psfrag,graphicx}
\usepackage{amsmath,times,amssymb}
\usepackage{cite}
\usepackage{mathrsfs}

\label{newcom}
\newcommand{\nbri}[2]{{#1}^{+{#2}}} 
\newcommand{\onbr}[2]{{#1}^{\bar{#2}}} 
\newcommand{\feas}[1]{{#1}^f} 
\newcommand{\opnbri}[2]{{#1}^{+{#2}*}} 
\newcommand{\aonbr}[2]{{#1}^{\bar{#2},a}} 

\newcommand{\pred}[1]{\textbf{#1}} 
\newcommand{\X}{\mathcal{X}}
\newcommand{\U}{\mathcal{U}}

\newcommand{\Nr}{\mathcal{N}}

\newcommand{\rr}{{\mathbb R}}
\newcommand{\CE}[1]{\mathop{CE}\left( #1 \right)}
\newcommand{\diag}[1]{\mathop{diag}\left( #1 \right)}
\newcommand{\cvd}{\hfill $\Box$}

\title{\LARGE \bf A {Jacobi} algorithm for distributed model predictive control of dynamically coupled systems}

\author{Dang Doan, Tam\'{a}s Keviczky, Ion Necoara and Moritz Diehl
\thanks{Dang Doan and Tam\'{a}s Keviczky are with the Delft Center for Systems and Control, Delft University of
Technology, Delft, The Netherlands. {\texttt
\small{minhdang@doan.vn}, \texttt \small{t.keviczky@tudelft.nl}}}
\thanks{Ion Necoara and Moritz Diehl are with the Department of
Electrical Engineering, Katholieke Universiteit Leuven, Leuven,
Belgium. {\texttt \small{
\{ion.necoara,moritz.diehl\}@esat.kuleuven.be}}} }

\begin{document}
\maketitle
\pagestyle{empty}
\thispagestyle{empty}

\begin{abstract} \label{abstract}
In this paper we introduce an iterative Jacobi algorithm for solving
distributed model predictive control (DMPC) problems, with linear
coupled dynamics and convex coupled constraints. The algorithm
guarantees stability and persistent feasibility, and we provide a
localized procedure for constructing an initial feasible solution by
constraint tightening. Moreover, we show that the solution of the
iterative process converges to the centralized MPC solution. The
proposed iterative approach involves solving local optimization
problems consisting of only few subsystems, depending on the choice
of the designer and the sparsity of dynamical and constraint
couplings. The gain in the overall computational load compared to
the centralized problem is balanced by the increased communication
requirements. This makes our approach more applicable to situations
where the number of subsystems is large, the coupling is sparse, and
local communication is relatively fast and cheap. A numerical
example illustrates the effects of the local problem size, and the
number of iterations on convergence to the centralized solution.
\end{abstract}

\section{Introduction} \label{intro}
Model predictive control (MPC) is the most successful advanced
control technology implemented in industry due to its ability to
handle complex systems with hard input and state constraints
\cite{Mac:02,MayRaw:00,GarPre:89}. The essence of MPC is to
determine a control profile that optimizes a cost criterion over a
prediction window and then to apply this control profile until new
process measurements become available. Then the whole procedure is
repeated and feedback is incorporated by using the measurements to
update the optimization problem for the next step.

For the control problem of large-scale networked systems,
centralized MPC is considered impractical, inflexible and unsuitable
due to information exchange requirements and computational aspects.
The subsystems in the network may belong to different authorities
that prevent sending all necessary information to one processing
center. Moreover, the optimization problem yielded by centralized
MPC can be excessively large for real-time computation. In order to
deal with these limitations, distributed MPC is proposed for control
of such large-scale systems, by decomposing the overall system into
small subsystems. The subsystems employ distinct MPC controllers,
use local information from neighboring subsystems, and collaborate
to achieve globally attractive solutions.

Approaches to distributed MPC design differ from each other in the
problem setup. In \cite{Camponogara:2002}, Camponogara \textit{et
al.} studied stability of coordination-based distributed MPC with
several information exchange conditions. In \cite{Dunbar:2006},
Dunbar and Murray proposed a distributed MPC scheme for problems
with coupled cost function, utilizing predicted trajectories of the
neighbors in each subsystem's optimization. Keviczky \textit{et al.}
proposed a distributed MPC scheme with a sufficient stability test
for dynamically decoupled systems in \cite{Keviczky:2006}, in which
each subsystem optimizes also over the behaviors of its neighbors.
Richards and How in \cite{Richards:2007} proposed a robust
distributed MPC method for networks with coupled constraints, based
on constraint tightening and a serial solution approach.

A distributed MPC scheme for dynamically coupled systems called
\emph{feasible-cooperation MPC} (FC-MPC) was proposed by Venkat
\textit{et al.} in \cite{Venkat:2005,Venkat:2005_techrep}, based on
a parallel synchronous approach for cooperative optimization. This
scheme works only for input-coupled linear time-invariant (LTI)
subsystem dynamics without state constraints, and is not applicable
to problems with constraints between subsystems. In this paper, we
propose an extension of this scheme in several ways in order to
solve these issues.

The distributed MPC algorithm described in this paper is able to
handle LTI dynamics with general dynamical couplings, and the
presence of convex coupled constraints. Each local controller
optimizes not only for itself, but also for its neighbors in order
to gain better overall performance. Global feasibility and stability
are achieved, whilst the algorithm can be implemented using local
communications. The proposed algorithm is based on an MPC framework
with zero terminal point constraint for increased clarity and
simplicity. We highlight an open research question that needs to be
addressed for a full treatment of the terminal cost based version of
this MPC framework, which would allow reduced conservativeness.
While other distributed MPC methods typically assume an initial
feasible solution to be available, we incorporate a decentralized
method to determine an initial feasible solution.

The problem formulation is described in Section~\ref{problem},
followed by two variations of the algorithm in
Section~\ref{algorithm}. It is shown that an algorithm using local
communication exists and it is equivalent to one that is based on
global communication. In Section~\ref{feasibility} 
we analyze the feasibility, stability and optimality of the
algorithm. Different ways of customizing the proposed algorithm and
a trade-off between communications and computational aspects are
discussed in Section~\ref{customizations}. Finally,
Section~\ref{experiment} illustrates the algorithm in a numerical
example and Section~\ref{conclusions} concludes the paper.

\section{Problem description} \label{problem}

\subsection{Coupled subsystem model}

Consider a plant consisting of $M$ subsystems. Each subsystem's
dynamics is assumed to be influenced directly by only a small number
of other subsystems. Let each subsystem be represented by a
discrete-time, linear time-invariant model of the form:
\begin{align} \label{eqn_dmodel}
x^i_{t+1} &= \sum_{j=1}^M (A_{ij} x^j_t + B_{ij} u^j_t), 
\end{align}
where $x^i_{t}\in \rr^{n_i}$ and $u^i_{t} \in \rr^{m_i}$ are the
states and control inputs of the $i$-th subsystem at time $t$,
respectively.

\begin{remark}
This is a very general model class for describing dynamical coupling
between subsystems and includes as a special case the combination of
\emph{decentralized models} and \emph{interaction models} in
\cite{Venkat:2005}.
\end{remark}

\begin{figure}
  \centering
  \includegraphics[width=.70\linewidth]{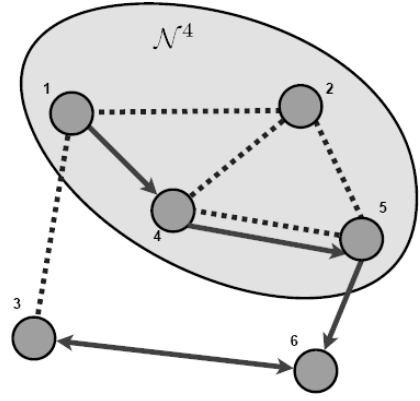}
  \caption[Interaction map]{An interaction map showing the constraints (dotted links) and dynamical couplings (arrows) between subsystems. In this example, $\Nr^4 = \{4, 1, 2 ,5 \}$. }
  \label{fig_Nr}
\end{figure}

We define the \emph{neighborhood of $i$}, denoted by $\Nr^i$, as the
set of indices of subsystems that have either direct dynamical or
convex constraint coupling with subsystem $i$. In
Figure~\ref{fig_Nr}, we demonstrate this with an \emph{interaction
map} where each node stands for one subsystem, the dotted links show
constraint couplings and the solid arrows represent dynamical
couplings. The neighborhood $\Nr^4$ of subsystem $4$ is the set of
$\{ 4, 1, 2, 5 \}$. We will refer to vectors and sets related to
nodes in $\Nr^i$ with a superscript $+i$. The collection of all
other nodes that are not included in $\Nr^i$ will be referred to
with a superscript $\bar{i}$.%

\subsection{Convex coupled constraints}

Each subsystem $i$ is assumed to have local convex coupled
constraints involving only a small number of the others. If we fix
the control inputs and the corresponding states of the nodes outside
$\Nr^i$, the state and input constraints involving the nodes in
$\Nr^i$ can be defined in the following way:
\begin{align} \label{eqn_constraints}
\nbri{x}{i}_{t} \in \nbri{\X}{i}(\onbr{x}{i}_t), \quad \nbri{u}{i}_t \in \nbri{\U}{i}(\onbr{u}{i}_t), \quad \forall i=1,\ldots, M
\end{align}
where $\nbri{\X}{i}(\onbr{x}{i}_t)$ and $\nbri{\U}{i}(\onbr{u}{i}_t)$ are
closed and convex sets parameterized by the states and control
inputs of nodes outside $\Nr^i$.
\begin{remark}
Note that the constraints involving nodes of $\Nr^i$ in general do
not depend on every other state and input outside $\Nr^i$, only on
the immediate neighbors of $\Nr^i$. The notation in
\eqref{eqn_constraints} is used for simplicity.
\end{remark}

\subsection{Centralized model}

Let $x = \left[ {x^1}^T \cdots {x^M}^T \right]^T$ and $u = \left[
{u^1}^T \cdots {u^M}^T \right]^T$ denote the aggregated states and
inputs of the full plant, with dimensions $\rr^{\sum_{i=1}^M n_i}$
and $\rr^{\sum_{i=1}^M m_i}$ respectively.
The matrices $A$ and $B$ will denote the aggregated subsystem
dynamics matrices and are assumed to be stabilizable:
\begin{align*}
A = \begin{bmatrix} A_{11} & \ldots & A_{1M} \\ \vdots & ~ & \vdots \\ A_{M1} & \ldots & A_{MM} \end{bmatrix},
B = \begin{bmatrix} B_{11} & \ldots & B_{1M} \\ \vdots &  & \vdots \\ B_{M1} & \ldots & B_{MM} \end{bmatrix}.
\end{align*}
The full (centralized) plant model is thus represented as:
\begin{align} \label{eqn_cmodel}
\begin{split}
x_{t+1} &= A x_t + B u_t, 
\end{split}
\end{align}


\begin{remark}
The \emph{centralized model} defined in \ref{eqn_cmodel} is more
general than the so-called \emph{composite model} employed in
\cite{Venkat:2005}. In our approach, the \emph{centralized model}
can represent both couplings in states and inputs. In
\cite{Venkat:2005}, the authors use an input-coupled \emph{composite
model}, which requires the subsystems' states to be decoupled,
allowing only couplings in inputs.
\end{remark}

\subsection{Centralized MPC problem}

The centralized MPC problem is formulated based on a typical
quadratic MPC framework \cite{Mac:02} with prediction horizon $N$,
and the following quadratic cost function at time step $t$:
\begin{align}
V_t = \sum_{k=0}^{N-1} {{x_{k,t}}^T Q x_{k,t} + {u_{k,t}}^T R u_{k,t}}
\label{eq_centralized_cost}
\end{align}
where $x_{k,t}$ denotes the centralized state vector at time $t + k$
obtained by starting from the state $x_{0,t} = x_t$ and applying to
system (\ref{eqn_cmodel}) the input sequence $u_{0,t}, \ldots,
u_{k-1,t}$. $Q=\diag{Q_1, \cdots, Q_M}$, $R=\diag{R_1, \cdots, R_M}$
with $\diag{.}$ function representing the block diagonal matrix.
Matrices $Q_i$ are positive semidefinite and $R_i$ are positive
definite.

Let $\pred{x}_t = [x_{1,t}^T, \cdots, x_{N,t}^T]^T$, $\pred{u}_t =
[u_{0,t}^T, \cdots, u_{N-1,t}^T]^T$.  The centralized MPC problem is
then defined as:
\begin{align}
  V_t^*(x_t) = \min_{\pred{x}_t, \pred{u}_t} ~&~ \sum_{k=0}^{N-1} x_{k,t}^T Q x_{k,t} + u_{k,t}^T R u_{k,t} \label{eq:cenMPC} \\
  \text{s.t.} ~&~ x_{k+1,t} = A x_{k,t} + B u_{k,t}, k=0,...,N-1, \nonumber \\
 &~ u_{k,t}  \in \U, k=0,...,N-1, \nonumber \\
 &~ x_{k,t}  \in \X, k=1,...,N-1, \nonumber \\
 &~ x_{N,t} = 0, \nonumber \\
 &~ x_{0,t} = x_t, \nonumber
\end{align}
where $\U$ and $\X$ are defined as $\bigcap_{i=1}^M \CE{\U^{+i}}$
and $\bigcap_{i=1}^M \CE{\X^{+i}}$, respectively. The $\CE{\cdot}$
operator denotes cylindrical extension to the set of
$\rr^{\sum_{i=1}^M n_i}$ and $\rr^{\sum_{i=1}^M m_i}$, respectively.
In other words, if $\X^{+i} \subset \rr^d_i$ then $\CE{\X^{+i}} =
\X^{+i} \times \rr^{\sum_{i=1}^M n_i - d_i}$. The vector $x_t$
contains the measured states at time step $t$.

Let $\pred{u}^*_t = [(u^*_{0,t})^T, \cdots, (u^*_{N-1,t})^T]^T$
denote the optimal control solution of (\ref{eq:cenMPC}) at time
$t$. Then, the first sample of $\pred{u}^*_t$ is applied to the
overall system:
\begin{align} \label{eq:mpclaw}
    u_t=u^*_{0,t}.
\end{align}
The optimization~(\ref{eq:cenMPC}) is repeated at time $t+1$, based
on the new state $x_{t+1}$. In \cite{KeeGil:88} it was shown that
with prediction horizon $N$ long enough to allow a feasible solution
to the optimization problem, the closed-loop system
(\ref{eqn_cmodel})-(\ref{eq:mpclaw}) is stable.

Before formulating the distributed MPC problems, we eliminate the
state variables in the centralized MPC formulation. In the following
we will also assume $t=0$ without loss of generality and drop
subscript $t$ for simplicity. The set of dynamics equations allows
us to write the predicted states as
\begin{align}
\pred{x} = \alpha \pred{u} + \beta(x_0),
\end{align}
where
\begin{align*}
\alpha = \begin{bmatrix}
B  & 0  & \ldots & 0 \\
AB & B  & \ddots & \vdots \\
\vdots     & \vdots & \ddots & 0 \\
A^{N-1}B & A^{N-2}B & \ldots & B
\end{bmatrix}, \quad
\beta (x_0) = \begin{bmatrix} A \\ A^2 \\ \vdots \\ A^N \end{bmatrix} x_0.
\end{align*}

Using the above equations, we can eliminate state variables in the
centralized MPC leading to the following problem:
\begin{align} \label{eqn_cmpc}
  \min_{\pred{u}} ~&~ V(\pred{u}, x_0) \\
  = \min_{\pred{u}} ~&~ \pred{u}^T ( \alpha^T \pred{Q} \alpha + \pred{R})\pred{u} + 2 (\alpha^T \pred{Q} \beta )^T \pred{u} + \beta^T \pred{Q} \beta \nonumber \\
  \text{s.t.} ~&~ \pred{u} \in \tilde{\U}, \nonumber \\
  &~ \alpha \pred{u} + \beta (x_0) \in \tilde{\X}, \nonumber \\
  &~ F \pred{u} + A^{N} x_0 = 0, \nonumber 
\end{align}
where $\tilde{\U} = \prod_{k=1}^N \U$ and $\tilde{\X} =
\prod_{k=1}^N \X$. $F=[A^{N-1}B, A^{N-2}B, \ldots, B]$ is the last
block row of $\alpha$. The matrices $\pred{Q}$ and $\pred{R}$ are
block-diagonal, built from weighting matrices $Q$ and $R$. Therefore
$\pred{Q}$ is positive semidefinite and $\pred{R}$ is positive
definite, making the cost function $V(\pred{u}, x_0)$ strictly
convex with any given $x_0$.

\subsection{Distributed MPC problem} \label{dmpc}

We will solve problem (\ref{eqn_cmpc}) by dividing it into smaller,
overlapping DMPC problems, with each DMPC assigned to one subsystem
but optimizing also over neighboring subsystems at the same time. At
each time step, by solving DMPCs and combining the local solutions
in an iterative process, we will get an increasingly accurate
approximate solution of the centralized MPC problem.

In the DMPC problem for subsystem $i$, the global cost function is
optimized with respect to a reduced set of variables: control inputs
of $i$ and its neighbors, denoted together by $\nbri{\pred{u}}{i}$. Each
DMPC problem will guarantee that all constraints of the centralized
MPC problem are satisfied. The DMPC of subsystem $i$ can be recast
into the following optimization problem:
\begin{align} \label{eqn_dmpc}
  \min_{\nbri{\pred{u}}{i}} ~&~ V (\pred{u}, x_0) \\
  \text{s.t.} ~&~ \pred{u} \in \tilde{\U}, \nonumber \\
  &~ \alpha \pred{u} + \beta (x_0) \in \tilde{\X}, \nonumber \\
  &~ F \pred{u} + A^{N} x_0 = 0, \nonumber \\
  &~ \onbr{\pred{u}}{i} = \aonbr{\pred{u}}{i} \nonumber
\end{align}
where $\aonbr{\pred{u}}{i}$ denotes the \emph{assumed inputs} of all
non-neighbors of $i$. For now we assume that in the beginning of
each step, each node $j$ transmits its \emph{assumed inputs} for
$\pred{u}^j = [(u^j_0)^T,...,(u^j_{N-1})^T]^T$ to the entire
network, node $i$ receive these vectors from $\forall j \not\in
\Nr^i$ to form $\aonbr{\pred{u}}{i}$.

Note that $\pred{u}$ is the combination of $\nbri{\pred{u}}{i}$ and
$\onbr{\pred{u}}{i}$. With each $i$, we can construct two pairs of
matrices $\nbri{\alpha}{i}$, $\onbr{\alpha}{i}$ and $\nbri{F}{i}$,
$\onbr{F}{i}$ so that:
\begin{align}
\alpha \pred{u} &= \nbri{\alpha}{i} \nbri{\pred{u}}{i} + \onbr{\alpha}{i} \onbr{\pred{u}}{i} \\
F \pred{u} &= \nbri{F}{i} \nbri{\pred{u}}{i} + \onbr{F}{i} \onbr{\pred{u}}{i} \nonumber
\end{align}

By eliminating the input constraints which involves only
$\onbr{\pred{u}}{i}$, the DMPC problem~(\ref{eqn_dmpc}) is
equivalent to the following
\begin{align} \label{eqn_dmpci}
  \min_{\nbri{\pred{u}}{i}} ~&~ V (\pred{u}, x_0) \\
  \text{s.t.} ~&~ \nbri{\pred{u}}{i} \in \nbri{\tilde{\U}}{i}, \nonumber \\
  &~ \nbri{\alpha}{i} \nbri{\pred{u}}{i} + \onbr{\alpha}{i} \onbr{\pred{u}}{i} + \beta (x_0) \in \tilde{\X}, \nonumber \\
  &~ \nbri{F}{i} \nbri{\pred{u}}{i} + \onbr{F}{i} \onbr{\pred{u}}{i} + A^{N} x_0 = 0, \nonumber \\
  &~ \onbr{\pred{u}}{i} = \aonbr{\pred{u}}{i} \nonumber
\end{align}
in which $\nbri{\tilde{\U}}{i} = \prod_{k=1}^N
\nbri{\U}{i}(\onbr{u}{i}_k)$. The optimal solution of
(\ref{eqn_dmpc}), (\ref{eqn_dmpci}) will be denoted by
$\opnbri{\pred{u}}{i}$.

For implementation, we introduce the notion of the \emph{$r$-step
extended neighborhood} for each subsystem $i$, denoted by $\Nr^i_r$,
which contains all nodes that can be reached from node $i$ in not
more than $r$ links. $\Nr^i_r$ is the union of subsystem indices in
the neighborhoods of all subsystems in $\Nr^i_{r-1}$ :
\begin{align}
\Nr^i_r = \bigcup_{j \in \Nr^i_{r-1} } \Nr^j,
\end{align}
where $\Nr^i_1 := \Nr^i$.
We see that in order to solve (\ref{eqn_dmpci}), subsystem $i$ only
needs information from other subsystems in $\Nr^i_{N+1}$, the
initial states and assumed inputs of subsystems outside
$\Nr^i_{N+1}$ are redundant.

\begin{remark}
The MPC formulation using terminal point constraint described above
simplifies our exposition but it is rather conservative. This could
be alleviated by using a dual-mode MPC formulation with a terminal
cost function. However, in order for this to be a truly distributed
approach, the terminal cost function associated with the terminal
controllers should have a sparse structure. This would allow the
construction of a centralized Lyapunov function in a local way,
using only local information. In \cite{Venkat:2005}, the authors try
to bypass this obstacle by using additional restrictive assumptions:
they employ zero terminal controllers and require all subsystems and
interaction models (coupled via the inputs only) to be stable. These
assumptions can actually be more conservative than using a terminal
point constraint, preventing the application of the FC-MPC method in
general dynamically coupled systems. Finding terminal controllers
that lead to a structured terminal cost is an open problem and
subject of our current research.
\end{remark}

\section{Jacobi-type algorithm} \label{algorithm}

In this section we present an iterative procedure to approximate the
centralized MPC solution by repeatedly calculating and combining the
solutions the local DMPC problems described in the previous section.
We will show two versions of our approach, which are based on Jacobi
distributed optimization \cite{BertTsit:parallel-comp-book}. The
proposed algorithms maintain feasibility of intermediate solutions
and converge to the centralized MPC solution asymptotically. The
first version uses global communication and can be considered as an
extension of FC-MPC \cite{Venkat:2005}. The second version relies on
local communication and represents the main contribution of the
paper. We will show that the two versions are equivalent to each
other, which leads to simplified analysis in
Section~\ref{feasibility}.

\subsection{Globally and locally communicating algorithms}

For each time step $t$, we assume that a feasible input
$\feas{\pred{u}}$ is given for the entire system.
(Section~\ref{subsec:localfeas} discusses a method of obtaining such
a feasible initial control sequence in a distributed way, given a
known initial condition.) In each step of the proposed DMPC scheme,
the subsystems cooperate and perform a Jacobi algorithm, where each
subsystem iteratively solves the optimization problem
(\ref{eqn_dmpci}) with regards to its local variables, and
incorporates a convex combination of
neighboring local solutions.  

During every MPC sampling period, a distributed iterative loop is
employed, and is indexed by $p$. At each iteration $p$,
$\feas{\pred{u}}$ is updated. We will refer to vectors obtained in
these iterations with subscript $(p)$.

For $p=1$, we initialize the iteration with $\pred{u}_{(0)} = \feas{\pred{u}}$. Let
$\pred{u}^{s|i}_{(p)}$ denote the control sequence of the whole system
stored in the memory of subsystem $i$ at iteration $p$. For making convex
combinations, each subsystem $i$ is assigned a weight $\lambda^i \in
(0,1)$, satisfying $\sum_{i=1}^{M} \lambda^i = 1$. The choice of
weights is arbitrary and could depend on the specific problem, the
simplest choice will be equal weights $(\lambda^i = \frac{1}{M},
\forall i)$. We propose then the following iterative algorithm:
\begin{algorithm}[Jacobi DMPC with global communication] \label{alg_dmpc}%

Given $N$, $p_{max} > 0$, $\epsilon >0$:

\begin{enumerate}
\item[\textbf{1.}] $p \leftarrow 1$, $\pred{u}_{(0)} \leftarrow \feas{\pred{u}},
\rho^i \leftarrow$ large number, $\forall i = 1, \ldots, M$.

\textbf{while} $\rho^i > \epsilon$ for some $i$ and $p \leq
p_{max}$

\begin{enumerate}
\item[\textbf{a.}] \textbf{for} each $i = 1,\ldots, M$ \\
Construct new $\aonbr{\pred{u}}{i}$ from $\pred{u}_{(p-1)}$.\\
Solve problem (\ref{eqn_dmpc}) to get
$\opnbri{\pred{u}}{i}_{(p)}$. Construct a global input vector
$\pred{u}^{s|i}_{(p)}$ from $\opnbri{\pred{u}}{i}_{(p)}$ and $\aonbr{\pred{u}}{i}$.
Transmit $\pred{u}^{s|i}_{(p)}$ to a central update location.%

\textbf{end (for)}
\item[\textbf{b.}] Merge local solutions according to the
    following convex combination:
\begin{align}
\pred{u}_{(p)} = \sum_{i=1}^M \lambda^i \; \pred{u}^{s|i}_{(p)}
\end{align}
\item[\textbf{c.}] Compute the progress and iterate:
\begin{align*}
\rho^i &= \| \pred{u}_{(p)} - \pred{u}_{(p-1)} \|
\end{align*}
\end{enumerate}
$p \leftarrow p+1$ \\
\textbf{end (while)}

\item[\textbf{2.}] Each subsystem implements the first input value in $\pred{u}^i_{(p)}$.
    \begin{align} \label{eq:dmpclaw}
        u^i_t=u^i_{0,(p)}.
    \end{align}

\item[\textbf{3.}] Shift the predicted input sequence one
    step to make a feasible solution for the following MPC
    update:
\begin{align*}
\feas{\pred{u}} = [u_{1,(p)}, \cdots, u_{N-1,(p)}, 0].
\end{align*}
\item[\textbf{4.}] $t \leftarrow t+1$. Measure new initial states $x_t$, go to step \textbf{1.}
\end{enumerate}
\end{algorithm}

Algorithm~\ref{alg_dmpc} requires the existence of a central
coordinator that communicates with all subsystems and performs the
convex combination to find $\pred{u}_{(p)}$. For implementation, a
control scheme without global communication is desired. Next we
introduce a variation of Algorithm~\ref{alg_dmpc} that only needs
local communication.

Let $\pred{u}^{i,f}$ denote the feasible input sequence of subsystem
$i$, and $\pred{u}^{j|i*}_{(p)}$ denote control sequence of
subsystem $j$ computed by subsystem $i$ when solving its DMPC
problem (\ref{eqn_dmpci}) at iteration $p$.

\begin{algorithm}[Jacobi DMPC with local communication\footnote{'local
communication' in this context means that each subsystem only
communicates with 'nearby' subsystems in a small region.}]
\label{alg_dmpci}

Given $N$, $p_{max} > 0$, $\epsilon >0$, and assuming each subsystem
$i$ knows a feasible input $\pred{u}^{j,f}$ for all subsystems $j
\in \Nr^i_{N+1}$.

\textbf{1.} $p \leftarrow 1$, $\rho^i \leftarrow$ large number, $\forall i = 1,
\ldots, M$.

\textbf{while} $\rho^i > \epsilon$ for some $i$ and $p \leq
p_{max}$

\textbf{for} each $i = 1,\ldots, M$
\begin{enumerate}
\item[\textbf{a.}] Subsystem $i$ solves the local problem
    (\ref{eqn_dmpci}), using $\{ \pred{u}^{j,f} | \forall j
    \in \Nr^i_{N+1} \backslash \Nr^i \}$ as assumed inputs
    for subsystems outside $\Nr^i$ but inside $\Nr^i_{N+1}$.
    The solution is comprised of $\{ \pred{u}^{j|i*}_{(p)}, \; j \in \Nr^i \}$.
\item[\textbf{b.}] Subsystem $i$ receives solutions for
    itself calculated by its neighbors $\{
    \pred{u}^{i|j*}_{(p)}, \; j \in \Nr^i \}$, then updates
    its solution for iterate $p$ according to:
\begin{align}
\pred{u}^{i}_{(p)} = \sum_{j \in \Nr^i} \lambda^j \pred{u}^{i|j*}_{(p)} + \left(1 - \sum_{j \in \Nr^i} \lambda^j\right) \pred{u}^{i}_{(p-1)}
\end{align}

\item[\textbf{c.}] Calculate the progress:
\begin{align*}
\rho^i &= \| \pred{u}^{i}_{(p)} - \pred{u}^{i}_{(p-1)} \|
\end{align*}

\item[\textbf{d.}] $\pred{u}^{i,f} \leftarrow
    \pred{u}^{i}_{(p)}$, subsystem $i$ transmits new
    $\pred{u}^{i,f}$ to all subsystems in $\Nr^i_{N+1}$.
\end{enumerate}

\textbf{end (for)}

$p \leftarrow p+1 $

\textbf{end (while)}

\textbf{2.} Each subsystem $i$ implements the first input
    value in $\pred{u}^i_{(p)}$:
    \begin{align} \label{eq:dmpclaw2}
        u^i_t=\pred{u}^i_{0,(p)}.
    \end{align}

\textbf{3.} Shift the predicted input sequence by one
    step to make a feasible solution for the following MPC
    update:
\begin{align*}
\pred{u}^{i,f} = [u^i_{1,(p)}, \cdots, u^i_{N-1,(p)}, 0], \;\; i = 1,\ldots, M .
\end{align*}

\textbf{4.} $t \leftarrow t+1 $. Measure new initial states $x^i_t$, go to step \textbf{1.}
\end{algorithm}

The major difference between Algorithms~\ref{alg_dmpc} and
\ref{alg_dmpci} is at step \textbf{1.b}: in Algorithm~\ref{alg_dmpc}
the convex combination is performed on the global control input
vector, while in Algorithm~\ref{alg_dmpci} each local controller
performs convex combination using its local control input vectors,
therefore removing the need of a coordinator. In the sequel, we will
show that the two algorithms are equivalent, thus allowing us to
implement Algorithm~\ref{alg_dmpci} while using
Algorithm~\ref{alg_dmpc} for analysis.

\subsection{Equivalence of the two algorithms}

The two crucial differences between Algorithm~\ref{alg_dmpci} and
\ref{alg_dmpc} are the communication requirement and the update
method. We already mentioned that the optimization problem
(\ref{eqn_dmpci}) is equivalent to (\ref{eqn_dmpc}), thus each
subsystem only has to transmit its new results to the subsystems
inside $\Nr^i_{N+1}$. This leads to the local communications in
Algorithm~\ref{alg_dmpci}. Now we will show that the local update of
Algorithm~\ref{alg_dmpci} is also equivalent to the global update of
Algorithm~\ref{alg_dmpc}.

Consider Algorithm~\ref{alg_dmpci}, at the beginning of iteration
$p$, a local input vector $\pred{u}^i_{p-1}$ is given for each $i$. Then each subsystem
$j \in \Nr^i$ computes $\pred{u}^{i|j*}_p$ and sends these solutions
to $i$, which forms the final update for itself $\pred{u}^{i}_p$. Note that $i \in \Nr^j
\Leftrightarrow j \in \Nr^i$, so we have
\begin{align} \label{eqn_localupdate}
\pred{u}^{i}_{(p)} &= \sum_{j \in \Nr^i} \lambda^j \pred{u}^{i|j*}_{(p)} + \left(1 - \sum_{j \in \Nr^i} \lambda^j\right) \pred{u}^{i}_{(p-1)} \nonumber \\
&= \sum_{j \in \Nr^i} \lambda^j \; \pred{u}^{i|j*}_{(p)}  - \left( \sum_{j \in \Nr^i} \lambda^j \right) \pred{u}^{i}_{(p-1)} + \pred{u}^{i}_{(p-1)}
\end{align}

Now consider Algorithm~\ref{alg_dmpc}, at the beginning of iteration
$p$, starting from the same $\pred{u}^i_{(p-1)}$ as in
Algorithm~\ref{alg_dmpci}. The local problem of the $j$-th
subsystem achieves solutions $\{ \pred{u}^{i|j*}_{(p)} ~|~ j \in \Nr^i
\}$ which are equal in both algorithms and are transmitted to
subsystem $i$. Then the global update $\pred{u}_{(p)}$ yields
$\pred{u}^{i}_{(p)}$ as follows:
\begin{align} \label{eqn_globalupdate}
\pred{u}_{(p)} &= \sum_{i=1}^M \lambda^i \; \pred{u}^{s|i}_{(p)} \nonumber \\
\Leftrightarrow \pred{u}^{i}_{(p)} &= \sum_{j \in \Nr^i} \lambda^j \; \pred{u}^{i|j*}_{(p)} + \sum_{j \not\in \Nr^i} \lambda^j \; \pred{u}^{i|j*}_{(p)} \nonumber \\
 &= \sum_{j \in \Nr^i} \lambda^j \; \pred{u}^{i|j*}_{(p)} + \sum_{j \not\in \Nr^i} \lambda^j \; \pred{u}^{i}_{(p-1)} \nonumber \\
 &= \sum_{j \in \Nr^i} \lambda^j \left( \pred{u}^{i|j*}_{(p)} - \pred{u}^{i}_{(p-1)} \right) + \sum_{j=1}^M \lambda^j \; \pred{u}^{i}_{(p-1)} \nonumber \\
 &= \sum_{j \in \Nr^i} \lambda^j \; \pred{u}^{i|j*}_{(p)} - \left( \sum_{j \in \Nr^i} \lambda^j \right) \pred{u}^{i}_{(p-1)} + \pred{u}^{i}_{(p-1)}
\end{align}

Comparing (\ref{eqn_localupdate}) and (\ref{eqn_globalupdate}) we
can see that the local update of Algorithm~\ref{alg_dmpci} and the
global update of Algorithm~\ref{alg_dmpc} yield the same result.
This implies that Algorithm~\ref{alg_dmpci} does exactly what
Algorithm~\ref{alg_dmpc} does, except it only needs to use regional
information (each subsystem $i$ needs to communicate with subsystems
in the region $\Nr^i_{N+1}$). If $M \gg N$ and the interaction map
is relatively sparse, this region will be much smaller than the
whole network, thus DMPC problems can be considered local
optimization problems.

The equivalence of the two proposed DMPC algorithms allows us to
prove their feasibility, stability and optimality aspects by
analyzing the globally communicating algorithm, which is more
comprehensive than the locally communicating algorithm. In the next
sections, we refer to Algorithm~\ref{alg_dmpc} for analysis.

\section{Feasibility, stability and optimality} \label{feasibility}

\subsection{Constructing initial feasible solutions locally} \label{subsec:localfeas}

Although in current literature it is typically assumed that an
initial centralized feasible solution exists and is available, in
this section we give a simple but implementable way of actually
constructing it in a distributed way assuming that the global
initial state is available in advance.

The initial feasible prediction input $\feas{\pred{u}}$ at time
$t=0$ can be calculated locally by using an inner approximation of
the global feasible set, which we will denote with $\Omega$ based on
all the constraints appearing in (\ref{eqn_cmpc}) and the global
initial state, which is assumed to be available. Consider an
inner-hyperbox (or hyperrectangular) approximation $\mathcal{B}$ of
the feasible set $\Omega$, which then takes the form of a Cartesian
product:
\begin{align} \label{eq:boxes}
\mathcal{B} = \mathcal{B}^1 \times \cdots \times \mathcal{B}^M, \quad \mathcal{B} \subset \Omega.
\end{align}

This approximation essentially decomposes and decouples the
constraints among subsystems by performing constraint tightening.
Each subsystem $i$ will thus have to include $\mathcal{B}^i$ in
their local problem setup. Since the Cartesian product of these
local constraint sets are included in the globally feasible set
$\Omega$, any combination of local solutions within
$\mathcal{B}^i$ will be globally feasible as well. Needless to say
that the local constraint sets that arise from this inner-hyperbox
approximation will be in general quite conservative, but at the same
time will allow the construction of a centralized feasible solution
locally to initialize Algorithm~\ref{alg_dmpc}.

Calculation of the inner-hyperbox approximation can be performed a
priori and the local $\mathcal{B}^i$ constraints distributed to each
subsystem. A polynomial-time procedure to compute a maximum volume
inner box of $\Omega$ could follow the procedure described in
\cite{BemFil:04}. Let us denote the dimension of the global input
vector with $d=\sum_{i=1}^M m_i$. If we represent a box as
$\mathcal{B}(\underline{u},\bar{u})=\{ u \in \rr^d ~|~ \underline{u}
\leq u \leq\bar{u} \}$, then
$\mathcal{B}(\underline{u}^*,\underline{u}^* + v^*)$ is a maximum
volume inner box of the full-dimensional polytope defined as
$\Omega=\{u \in \rr^d ~|~ Au \leq b \}$, where
$(\underline{u}^*,v^*)$ is an optimal solution of
\begin{align} \label{eq:innerbox}
\begin{split}
\max_{\underline{u}, v} ~&~ \sum_{j=1}^{d} \ln v_j \\
\text{s.t.}~&~A\underline{u} + A^{+}v \leq b.
\end{split}
\end{align}
The matrix $A^{+}$ is the positive part of $A$. Obtaining the local
component-wise constraints $\mathcal{B}^i$ is then straightforward.
For time steps other than $t=0$, we construct a feasible solution by
performing step \textbf{3} of Algorithm~\ref{alg_dmpc}.

\subsection{Maintaining feasibility throughout the iterations}

Observe that in step \textbf{1.a}, we get $M$ feasible solutions
$\pred{u}^{s|i}_{(p)}$ for the centralized problem (\ref{eqn_cmpc}).
In step \textbf{1.b}, we construct the new control profile
$\pred{u}_{(p)}$ as a convex combination of these solutions. Since problem (\ref{eqn_cmpc}) is a convex
constrained QP, any convex combination of $\{ \pred{u}^{s|i}_{(p)} \}_{i=1}^{M}$
also satisfies the convex constraint set. Therefore $\pred{u}_{(p)}$ is
a feasible solution of optimization problems (\ref{eqn_cmpc}), and (\ref{eqn_dmpc}) for all $i$.

\subsection{Stability analysis} \label{stability}

Showing stability of the closed-loop system
(\ref{eqn_cmodel})-(\ref{eq:dmpclaw}) follows standard arguments for
the most part \cite{KeeGil:88,MayRaw:00}. In the following, we
describe only the most important part for brevity, which considers
the nonincreasing property of the value function. The proof in this
section is closely related to the stability proof of the FC-MPC
method in \cite{Venkat:2005_techrep}, the main difference is due to
the underlying MPC schemes: this method uses terminal point
constraint MPC while FC-MPC uses dual-mode MPC.

Let $\bar{p}_t$ and $\bar{p}_{t+1}$ stand for the last iteration
number of Algorithm~\ref{alg_dmpc} at step $t$ and $t+1$,
respectively. Let $V_{t}=V (\pred{u}_{(\bar{p}_t)}, x_t)$ and
$V_{t+1}=V (\pred{u}_{(\bar{p}_{t+1})}, x_{t+1})$ denote the cost
values associated with the final combined solution at step $t$ and
$t+1$, respectively. At step $t+1$, let $\Phi^{i}_{(p+1)} = V
(\pred{u}^{s|i}_{(p+1)}, x_{t})$ denote the global cost associated
with solution of subsystem $i$ at iterate $p+1$, and $\Phi_{(p)} = V
(\pred{u}_{(p)}, x_{t})$ the cost corresponding to the combined
solution at iterate $p$.

The global cost function can be used as a Lyapunov function, and its
nonincreasing property can be shown following the chain:
\begin{align*}
V_{t+1} \leq \cdots \leq
\Phi_{(p+1)} &\leq \Phi_{(p)} \leq \cdots \\
\cdots &\leq \Phi_{(1)}
\leq V_{t} - x_{t}^T Q x_{t} - u_{t}^T R u_{t}
\end{align*}

The two main components of the above inequality chain are shown in
the following two subsections.

\textbf{Showing that $\Phi_{(p+1)} \leq \Phi_{(p)}$}

The cost $V (\pred{u}, x_t)$ is a convex function of $\pred{u}$, thus
\begin{align}
\label{ineq_Phi_p_1}
V \left( \sum_{i=1}^M \lambda^i \; \pred{u}^{s|i}_{(p+1)}, x_t \right) \leq \sum_{i=1}^M \lambda^i \; V \left( \pred{u}^{s|i}_{(p+1)}, x_t \right)
\end{align}
Moreover, each $\pred{u}^{s|i}_{(p+1)}$ is the optimizer of $i$-th
local problem starting from $\pred{u}_{(p)}$, therefore we have:
\begin{align}
\label{ineq_Phi_i}
V \left( \pred{u}^{s|i}_{(p+1)}, x_t \right) = \Phi^i_{(p+1)} \leq \Phi_{(p)}, i = 1,\ldots, M
\end{align}
Substituting (\ref{ineq_Phi_i}) into (\ref{ineq_Phi_p_1}) leads to:
\begin{align*}
\Phi_{(p+1)} = V \left(\sum_{i=1}^M \lambda^i \; \pred{u}^{s|i}_{(p+1)}, x_t \right) \leq \sum_{i=1}^M \lambda^i \; \Phi_{(p)} = \Phi_{(p)}
\end{align*}
Using the above inequality, we can trace back to $p=1$:
\begin{align*}
V_{t+1} \leq \cdots. \leq \Phi_{(p+1)} \leq \Phi_{(p)} \leq \cdots \leq \Phi_{(1)}.
\end{align*}

\textbf{Showing that $\Phi_{(1)} \leq V_{t} - x_{t}^T Q x_{t} - u_{t}^T R u_{t}$}

At step $t+1$ and iteration $p=1$, recall that the initial feasible
solution $\feas{\pred{u}}$ of the centralized MPC is built by
Algorithm~\ref{alg_dmpc} at the end of step $t$ in the following
way:
\begin{align*}
\feas{\pred{u}} = [u_{1,(\bar{p}_t)}, \cdots, u_{N-1,(\bar{p}_t)}, 0]
\end{align*}

The DMPC of each subsystem $i$ optimizes the cost with respect to
$\nbri{\pred{u}}{i}$ starting from $\feas{\pred{u}}$, therefore
$\forall i = 1,\ldots, M$:
\begin{align*}
V \left( \pred{u}^{s|i}_{(1)}, x_t \right) & \leq V \left( \feas{\pred{u}}, x_t \right) \\
\Leftrightarrow \;\; \Phi^{i}_{(1)} & \leq \sum_{k=1}^{N-1} \left( x_{k,(\bar{p}_t)}^T Q x_{k,(\bar{p}_t)} + u_{k,(\bar{p}_t)}^T R u_{k,(\bar{p}_t)} \right) \\
\Leftrightarrow \;\; \Phi^{i}_{(1)} & \leq \Phi (\pred{u}_{(\bar{p}_t)}) - x_{0,\bar{p}_t}^T Q x_{0,(\bar{p}_t)} - u_{0,(\bar{p}_t)}^T R u_{0,(\bar{p}_t)} \\
\Leftrightarrow \;\; \Phi^{i}_{(1)} & \leq V_{t} - x_{t}^T Q x_{t} - u_{t}^T R u_{t}.
\end{align*}
%

Moreover, due to the convexity of $V (\pred{u}, x_t)$ and the convex
combination update $\pred{u}_{(1)} = \sum_{i=1}^M \lambda^i
\pred{u}^{s|i}_{(1)}$, we obtain
\begin{align*}
\Phi_{(1)} &= V \left( \sum_{i=1}^M \lambda^i \; \pred{u}^{s|i}_{(1)}, x_t \right) \leq \sum_{i=1}^M \lambda^i \; V \left( \pred{u}^{s|i}_{(1)}, x_t \right) \\
\Rightarrow \Phi_{(1)} & \leq \sum_{i=1}^M \lambda^i \; \Phi^{i}_{(1)} \leq \sum_{i=1}^M \lambda^i \; [V_{t} - x_{t}^T Q x_{t} - u_{t}^T R u_{t}] \\
\Leftrightarrow \Phi_{(1)} & \leq V_{t} - x_{t}^T Q x_{t} - u_{t}^T R u_{t}
\end{align*}

The above inequalities show that the value function decreases along
closed-loop trajectories of the system. The rest of the proof for
stability follows standard arguments found for instance in
\cite{Venkat:2005,KeeGil:88}.

\subsection{Optimality analysis} \label{optimality}

Using the descent approach, we will show that the solution of
Algorithm~\ref{alg_dmpc} approaches the solution of the centralized
MPC in (\ref{eqn_cmpc}), as $p \rightarrow \infty$. We characterize
the optimality of the proposed iterative procedure by using the
following results.

\begin{lemma}
A limit point of $\{ \pred{u}_{(p)} \}$ is guaranteed to exist.
\end{lemma}

\emph{Proof}: The feasible set of (\ref{eqn_dmpc}) is compact. It is
shown that every $\pred{u}_{(p)}$ is feasible, therefore this
sequence is bounded, thus converges. \cvd

\begin{lemma}
Every limit point of $\{ \pred{u}_{(p)} \}$ is an optimal solution of (\ref{eqn_cmpc}).
\end{lemma}

\textbf{Proof}: We will make use of the strict convexity of $V
(\cdot)$ and a technique, which is inspired by the proof of
Gauss-Seidel distributed optimization algorithms in
\cite{BertTsit:parallel-comp-book}. In our context however, we
address also the overlapping variables that are present in the local
optimization problems.

Let $\pred{v} = (\pred{v}^1,...,\pred{v}^M)$ be a limit point of $\{
\pred{u}_{(p)} \}$, assume $\{ \pred{u}_{(p')} \}$ is a subsequence
of $\{ \pred{u}_{(p)} \}$ that converges to $\pred{v}$.

In the following, we drop parameter $x_t$ in $V(\cdot)$ for
simplicity. Using the continuity of $V(\cdot)$ and the convergence
of $\{ \pred{u}_{(p')} \}$ to $\pred{v}$, we see that $V \left(
\pred{u}_{(p')} \right)$ converges to $V(\pred{v})$. This implies
the entire sequence $\{ V \left( \pred{u}_{(p)} \right) \}$
converges to $V(\pred{v})$. It now remains to show that $\pred{v}$
minimizes $V(\cdot)$ over the feasible set of (\ref{eqn_cmpc}).

We first show that $\opnbri{\pred{u}}{1}_{(p' + 1)} -
\nbri{\pred{u}}{1}_{(p')}$ converges to zero. Recall that at
iteration $p$, $\nbri{\pred{u}}{1}_{(p')}$ and
$\aonbr{\pred{u}}{1}_{(p')}$ forms $\pred{u}_{(p')}$, at iteration
$p+1$, $\opnbri{\pred{u}}{1}_{(p' + 1)}$ and
$\aonbr{\pred{u}}{1}_{(p')}$ forms $\pred{u}^{s|1}_{(p' + 1)}$.
Assume the contrary, or $\pred{u}^{s|1}_{(p' + 1)} -
\pred{u}_{(p')}$ does not converge to zero. There exists some
$\epsilon > 0$ such that $\| \pred{u}^{s|1}_{(p' + 1)} -
\pred{u}_{(p')} \| \geq \epsilon$ for all $p'$.

Let us fix some $\gamma \in (0,1)$ and define
\begin{align}
\pred{s}^1_{(p')} = \gamma \pred{u}_{(p')} + (1-\gamma)
\pred{u}^{s|1}_{(p' + 1)},
\end{align}
this means $\pred{s}^1_{(p')}$ lies between $\pred{u}_{(p')}$ and
$\pred{u}^{s|1}_{(p' + 1)}$, only differs from them in the values of
$\nbri{\pred{u}}{1}$.

Notice that $\pred{s}^1_{(p')}$ belongs to a compact set and
therefore has a limit point, denoted $\pred{s}^1_{\infty}$. Since
$\gamma \neq 1$ and $\pred{u}_{(p')} \neq \pred{u}^{s|1}_{(p' + 1)},
\forall p'$, we have $\pred{s}^1_{\infty} \neq \pred{v}$.

Using convexity of $V(\cdot)$, we obtain
\begin{align}
V \left( \pred{s}^1_{(p')} \right) \leq \max \{ V \left(
\pred{u}_{(p')} \right), V \left( \pred{u}^{s|1}_{(p' + 1)} \right)
\}.
\end{align}

Be definition, $\pred{u}^{s|1}_{(p' + 1)}$ minimizes $V(\cdot)$ over
the subspace of $\nbri{\pred{u}}{1}$. So we have:
\begin{align} \label{eqn_bound_V}
V \left( \pred{u}^{s|1}_{(p' + 1)} \right) \leq V \left(
\pred{s}^1_{(p')} \right) \leq V \left( \pred{u}_{(p')} \right)
\end{align}.

From Section~\ref{stability}, we have
\begin{align}
V \left( \pred{u}^{s|i}_{(p'+1)} \right) &\leq V \left( \pred{u}_{(p')} \right), ~~ \forall i=0,\cdots,M \label{V_s_i_1} \\
V \left( \pred{u}_{(p'+1)} \right) &\leq \sum_{i=1}^M \lambda^i V
\left( \pred{u}^{s|i}_{(p'+1)} \right) \label{V_s_i_2}
\end{align}

Taking the limit of (\ref{V_s_i_1}) and (\ref{V_s_i_2}), we obtain
\begin{align}
\lim_{p' \to \infty} V \left( \pred{u}^{s|i}_{(p'+1)} \right) =
V(\pred{v}), ~~ \forall i=0,\cdots,M
\end{align}

As $V \left( \pred{u}_{(p')} \right)$ and $V \left(
\pred{u}^{s|1}_{(p' + 1)} \right)$ both converge to $V(\pred{v})$,
taking limit of (\ref{eqn_bound_V}), we conclude that $V(\pred{v}) =
V(\pred{s}^1_{\infty})$, for any $\gamma \in (0,1)$. This
contradicts the strict convexity of $V(\cdot)$ in the subspace of
$\nbri{\pred{u}}{1}$. The contradiction establishes that
$\opnbri{\pred{u}}{1}_{(p' + 1)} - \nbri{\pred{u}}{1}_{(p')}$
converges to zero, leading to the convergence of
$\pred{u}^{s|1}_{(p' + 1)}$ to $\pred{v}$.

We have, by definition
\begin{align}
V \left( \pred{u}^{s|1}_{(p' + 1)} \right) \leq V \left(
\nbri{\pred{u}}{1}, \onbr{\pred{u}}{1}_{(p')} \right)
\end{align}

Taking the limit as $p'$ tends to infinity, we obtain
\begin{align}
V (\pred{v}) \leq V \left( \nbri{\pred{u}}{1}, \onbr{\pred{v}}{1}
\right)
\end{align}
or $\pred{v}$ is optimizer of $V(\cdot)$ in the subspace of
$\nbri{\pred{u}}{1}$. If we further consider $V (\cdot)$ in a
subspace corresponding to $\pred{u}^1$, then $V (\pred{v})$ is still
a minimum. Thus, the necessary optimality condition gives
${\nabla_{\pred{u}^1} V (\pred{v})}^T ( \pred{u}^1 - \pred{v}^1 )
\geq 0, \forall \pred{u}^1 \in \Omega_1$ where $\Omega_1$ is the
feasible set of (\ref{eqn_dmpc}) with $i=1$.

Repeating the procedure, we obtain
\begin{align}
\label{partial_opt} {\nabla_{\pred{u}^i} V (\pred{v})}^T (
\pred{u}^i - \pred{v}^i ) \geq 0,
\end{align}
for all $\pred{u}^i$ such that $\pred{u}^i$ is a feasible solution
of (\ref{eqn_dmpc}).

By summing up the system of equations in (\ref{partial_opt}) for
$i=1, \cdots, M$, we get:
\begin{align}
{\nabla_{\pred{u}} V(\pred{v})}^T (\pred{u} - \pred{v}) \geq 0,
\end{align}
for all $\pred{u}$ that is a feasible solution of (\ref{eqn_cmpc}).

This shows that $\pred{v}$ satisfies the optimality condition of
problem (\ref{eqn_cmpc}).\cvd

Using strict convexity of $V (\cdot)$, it follows that $\pred{v}$ is
in fact the global optimizer of Algorithm~\ref{alg_dmpc}.

\section{Communications and computational aspects} \label{customizations}

In this section, we discuss the communications and computational
aspects of our approach and illustrate the freedom that the designer
has in choosing the appropriate trade-off and performance level in a
certain application.

Although the overall computational load is reduced by employing the
distributed Algorithm~\ref{alg_dmpci}, its iterative nature implies
that communication between neighboring systems increases in
exchange. This trade-off is illustrated in Table~\ref{tb_sm_comms},
which compares the communication requirements of the centralized and
our distributed MPC approach. This overview suggests that our scheme
is mostly applicable in situations where local communication is
relatively fast and cheap.

\begin{table}[htb]
\caption{Comparison of communications requirements} 
\begin{center}
\begin{tabular}{|c|c|c|}
\hline
  & Centralized MPC & Distributed MPC \\
\hline
Communication & Global & Local \\
\hline
Each subsystem &  & Other subsystems \\
communicates & Central coordinator &  in $(N+1)$-step \\
with & & extended neighborhood \\
\hline
Total number of & & \\
messages sent & $2 \times M $ & $p_{max} \times 2 \times \sum_{i=1}^{M}{| \Nr^i_{N+1} |}$ \\
in each time step & & \\
\hline
\end{tabular}
\end{center}
\label{tb_sm_comms}
\end{table}

\begin{table}[htb]
\caption{Comparison of optimization problems} 
\begin{center}
\begin{tabular}{|c|c|c|}
\hline
  & Centralized MPC & Distributed MPC \\
\hline
Number of variables in & & \\
one optimization problem& $N \times M$ &  $N \times | \Nr^i |$ \\
\hline
Number of optimizations & & \\
solved in one sampling period & 1 & $p_{max} \times M$ \\
\hline
\end{tabular}
\end{center}
\label{tb_sm_opt_size}
\end{table}

Table~\ref{tb_sm_opt_size} shows the difference in size of the
optimization problems solved by the distributed and the centralized
method. Since $| \Nr^i |\ll M $, where $M$ is the total number of
subsystems, the local optimization problems in DMPC are much smaller
than the centralized one. 
Note that during one sampling period, the local DMPC optimization
problems are solved at most $p_{max}$ times. Nevertheless, DMPC is
in general more computationally efficient than the centralized MPC,
with a proper choice of $p_{max}$.

Choosing an appropriately high $p_{max}$ value leads to better
performance of the whole system.
The trade-off is that the increase of $p_{max}$ will also lead to
increased communications, and more local optimization problems will
need to be solved in one sampling time.

Another way to customize the algorithm is to expand the size of the
neighborhood that each subsystem optimizes for. In the proposed
Algorithms~\ref{alg_dmpc},~\ref{alg_dmpci}, each subsystem optimizes
for its direct neighbors when solving the local optimization. We may
have better performance when each subsystem optimizes also for its
2, 3, k-step expanded neighbors. Although we do not provide a formal
proof of this, we will give an illustration in the following section
through a numerical example. The intuition behind this phenomenon is
nevertheless clear: each subsystem will have more precise
predictions when it takes into account the behaviors of more
neighboring subsystems.

\section{Numerical example} \label{experiment}

\begin{figure}
  \centering
  \includegraphics[width=.99\linewidth]{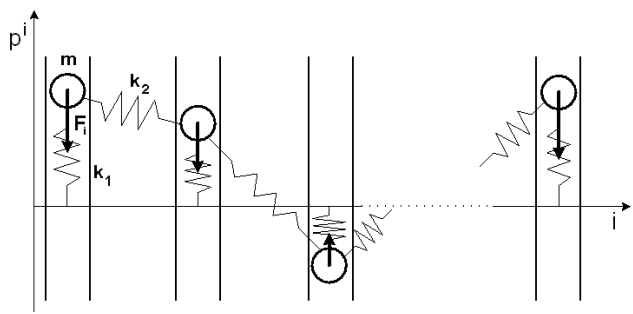}
   \caption{Setup of coupled oscillators}
   \label{fig_setup}
\end{figure}

In this section, we illustrate the application of
Algorithm~\ref{alg_dmpci} to a problem involving coupled
oscillators. The problem setup consists of $M$ oscillators that can
move only along the vertical axis, and are coupled by springs that
connect each oscillator with its two closest neighbors. An exogenous
vertical force will be used as the control input for each
oscillator. The setup is shown in Figure~\ref{fig_setup}.

Each oscillator is considered as one subsystem. Let the superscript
$i$ denote the index of oscillators. The dynamics equation of
oscillator $i$ is then defined as
\begin{equation} \label{eqn_dynamics}
m a^i = k_1 p^i - f_s v^i + k_2(p^{i-1} - p^i) + k_2(p^{i+1} - p^i) + F^i,
\end{equation}
where $p^i$, $v^i$ and $a^i$ denote the position, velocity and
acceleration of oscillator $i$, respectively. The control force
exerted at oscillator $i$ is $F^i$ and the parameters are defined as

$k_1$: stiffness of vertical spring at each oscillator

$k_2$: stiffness of springs that connect the oscillators

$m$: mass of each oscillator

$f_s$: friction coefficient of movements

From some nonzero initial state, the system needs to be stabilized
subject to the constraints:
\begin{equation}
\begin{vmatrix} p^i - \frac{p^{i-1} + p^{i+1}}{2} \end{vmatrix} \leq 4, \;\; \forall i=2,...,M-1
\end{equation}

Based on dynamical couplings and constraint couplings, the
neighborhood of each subsystem inside the chain is defined to
contain itself and its two closest neighbors $\Nr^i = \{ i-1, i, i+1
\}, i=2,...,M-1$, while for the two ends $\Nr^1 = \{ 1, 2 \}$ and
$\Nr^M = \{ M, M-1 \}$. We define the state vector as $x^i = [p^i,
v^i]^T$, and the input as $u^i = F^i$. The discretized dynamics with
sampling time $T_s$ is represented by the following matrices:
\begin{align*}
A_{ij} &= \left[ \begin{array}{cc}
0 & 0 \\
0 & 0 \end{array} \right], \forall j \not\in \Nr^i \\
A_{i,i-1} &= \left[ \begin{array}{cc}
0 & 0 \\
T_s k_2 & 0 \end{array} \right], \forall i=2,...,M \\
A_{ii} &= \left[ \begin{array}{cc}
1 & T_s \\
T_s (k_1 - 2 k_2) & 1 - T_s f_s \end{array} \right], \forall i=1,...,M \\
A_{i,i+1} &= \left[ \begin{array}{cc}
0 & 0 \\
T_s k_2 & 0 \end{array} \right], \forall i=1,...,M-1 \\
B_{ij} &= \left[ \begin{array}{c}
0 \\
0 \end{array} \right] \forall j \neq i \\
B_{ii} &= \left[ \begin{array}{c}
0 \\
T_s \end{array} \right], \forall i=1,...,M
\end{align*}

The following parameters were used in the simulation example:
\begin{align*}
k_1 &= 0.4, \quad k_2 = 0.3 \\
f_s &= 0.4, \quad T_s = 0.05, \quad m = 1 \\
M &= 40, \quad N = 20 \\
Q_i &= \left[ \begin{array}{cc}
100 & 0 \\
0 & 0 \end{array} \right], \quad R_i = 10
\end{align*}

\begin{figure}
  \centering
  \includegraphics[width=.90\linewidth]{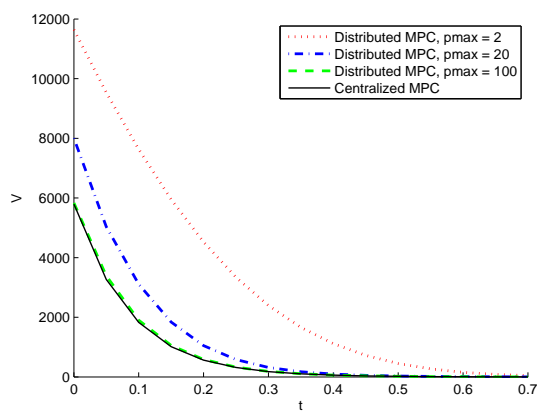}
  \caption[Comparison of global cost in centralized and distributed MPC simulations]
  {Time evolution of the global cost value of the centralized MPC in comparison with
  the distributed MPC algorithm using $p_{max}=2$, $p_{max}=20$ and $p_{max}=100$.}
  \label{fig_compare_cost}
\end{figure}

Starting from the same feasible initial state, we apply
Algorithm~\ref{alg_dmpci} with $p_{max}=2$, $20$ and $100$. The
results are compared to the solution obtained from the centralized
MPC approach. The results indicate that all states of the $40$
subsystems are stabilized. Figure~\ref{fig_compare_cost} shows the
evolution of the overall cost achieved by DMPC compared to the cost
of the centralized approach. We can see that the difference is
reduced by choosing a larger $p_{max}$ value.
Our analysis guarantees in fact that the DMPC solution converges to
the centralized one as $p$ tends to infinity.

\begin{figure}
  \centering
  \includegraphics[width=.90\linewidth]{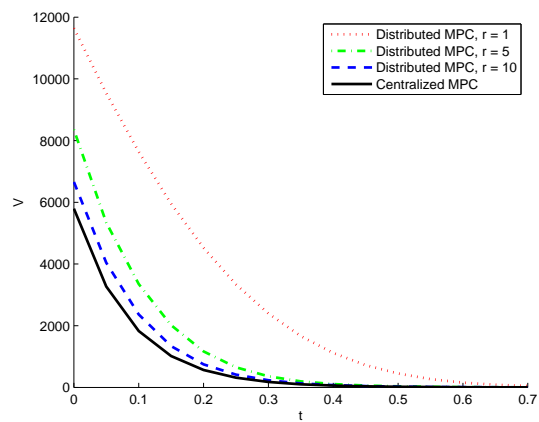}
  \caption{Time evolution of the global cost value of distributed MPC
  algorithms with different radius of neighborhood to be optimized by one local controller.}
  \label{fig_compare_r}
\end{figure}

As mentioned above, another way to customize the proposed
distributed MPC algorithm is for each local problem to consider
optimizing over the inputs of subsystems in a larger neighborhood.
Figure~\ref{fig_compare_r} illustrates the effect of optimizing in
each subproblem over an $r$-step extended neighborhood, with
$r=\{1,5,10\}$. Fixing the number of maximum subiterations to
$p_{max}=2$, we can observe a steady improvement in performance
until the increased neighborhood of each subsystem covers
essentially all other subsystems and end up with a centralized
problem.

\section{Conclusions} \label{conclusions}

We presented a Jacobi algorithm for solving distributed model
predictive control problems, which is able to deal with general
linear coupled dynamics and convex coupled constraints. We
incorporated neighboring subsystem models and constraints in the
formulation of the local problems for enhanced performance. Global
feasibility and stability were achieved, and a local implementation
of the algorithm was given, which relies on information exchange
from an extended set of ``nearby'' neighboring subsystems. It was
shown that the distributed MPC solution converges to the centralized
one through a localized iterative procedure. An a priori
approximation procedure was proposed, which allows to construct an
initial feasible solution locally by tightening constraints. We also
discussed the trade-off between communications and computational
aspects, the effect of increasing the maximum number of iterations
($p_{max}$) in one sampling period and the potential improvements
that can be gained by incorporating several subsystems into a local
optimization. We are currently working on an extension of the
algorithm, which allows the use of terminal costs in a dual-mode MPC
formulation.

\bibliographystyle{plain}
\bibliography{Dang_thesis}
\end{document}